\documentclass[11pt]{amsart}

\usepackage{eucal,amssymb,enumerate,enumitem}
\usepackage{graphics}

\input xy
\xyoption{all}
\newtheorem{thrm}{Theorem}[section]

\newtheorem{prp}[thrm]{Proposition}

\newtheorem{dfn}[thrm]{Definition}

\newtheorem{rmk}[thrm]{Remark}

\begin{document}

\title{A look at the prime and semiprime operations of one-dimensional domains}
\author{Janet C. Vassilev}
\address{The University of New Mexico, Albuquerque, NM 87131}
\email{jvassil@math.unm.edu} \subjclass[2000]{13A15, 13B22, 13C99}
\date{\today}

\keywords{closure operation, semiprime operation, prime operation,
star operation, integral closure}

\begin{abstract}{We continue the analysis of prime and semiprime
operations over one-dimensional domains started in \cite{Va}.  We first show that there are no bounded semiprime operations on the set of
fractional ideals of a one-dimensional domain.  We then prove the only prime operation is the identity on the set of ideals in semigroup
rings where the ideals are minimally generated by two or fewer elements.  This is not likely the case in semigroup rings with ideals of
three or more generators since we are able to exhibit that there is a non-identity prime operations on the set of ideals of
$k[[t^3,t^4,t^5]]$. }
 \end{abstract}

\maketitle

\section{Introduction}

In his 1935 book Idealtheorie \cite{Kr1}, Krull  defined an operation $I
\mapsto I^\prime$ on the set of fractional ideals of a domain $R$ to be a
{\it $\prime$-operation} if it satisfies the following properties where
$I$ and $J$ are ideals and $b$ is a regular element:
\begin{enumerate}
\item $I \subseteq I^\prime$
\item If $I \subseteq J$, then $I^\prime \subseteq J^\prime$.
\item $(I^\prime)^\prime=I^\prime$.
\item $I^\prime J^\prime \subseteq (IJ)^\prime$
\item $(bI)^\prime=bI^\prime$
\item $I^\prime +J^\prime \subseteq (I+J)^\prime$.
\end{enumerate}
Then a year later in \cite{Kr2}, he discussed the integral
completion or $b$-operation in terms of $\prime$-operations and
mentioned that he left out the properties:

\begin{enumerate}[resume]
\item  $R=R^\prime$
\item  $(I^\prime \cap J^\prime )^\prime =I^\prime \cap J^\prime $.
\end{enumerate}

In fact, Sakuma \cite{Sa} shows in 1957 that when looking at prime
operations on the set of fractional ideals of a domain, properties
(4), (6) and (8) are consequences of properties of (1), (2), (3),
(5) and (7). In 1964, Petro \cite{Pe} called the operations
satisfying properties (1)-(4) on the set of fractional ideals,
semiprime operations.  After, Gilmer wrote the book, Multiplicative
Ideal Theory \cite{Gi} in 1972, prime operations on the set of
fractional ideals came to be known as star operations.  Gilmer
explained in his book that $I^\prime$ was a common notation for
integral closure and as Krull also used the notation $I^\star$ with
regard to the $v$-operation, he preferred to to call these
operations star operations or $\star$-operations.  For more on star
operations and the structure on the set of star operations, one may
refer to the following two articles \cite{AA} by Dan and Dave
Anderson and \cite{AC} by Dan Anderson and Sylvia Cook.

In 1969, Kirby \cite{Ki} defined a closure operation on the set of
submodules of an $R$-module $M$ for an arbitrary commutative ring to
satisfy properties similar to (1), (2) and (3) above and
$(IN)^\prime=IN^\prime$ where $I \subseteq R$ is an ideal and $N$ is a
$R$-submodule of $M$. This seems to be the first reference where the such
operations were defined for an arbitrary commutative ring and on a set
other than the set of fractional ideals over a domain. The terms prime and
semiprime operation were reintroduced on the set of ideals of a
commutative ring by Ratliff in his 1989 paper \cite{Ra} on
$\Delta$-closures of ideals. Heinzer, Ratliff and Rush \cite{HRR} also use
the term semiprime operation when referring to the basically full closure
on the set of $m$-primary submodules of a module over a local ring
$(R,m)$.

Certainly, when the ring $R$ is a domain, we can determine which
semiprime operations on the set of ideals of $R$ are star operations
when we extend them to the set of fractional ideals. In a recent
paper \cite{Va}, the author determined all the semiprime operations
on the set of ideals of a Dedekind domain and all the semiprime
operations on the set of ideals of the complete cuspidal cubic. We
will observe that these operations are not semiprime operations on
the set of all fractional ideals; hence, they are not star
operations. Also, we exhibit that over $k[[t^3,t^4,t^5]]$, there is
a prime operation that is not the identity and it seems likely that
a one-dimensional semigroup ring having ideals minimally generated
by 3 or more elements, will have prime operations which are not the
identity. However, if the ideals are minimally generated by two or
fewer elements, then the only prime operation is the identity.

\section{There are no bounded semiprime operations on the set of fractional
ideals of a one-dimensional domain}

Let $R$ be an integral domain and $K$ its field of fractions.  Let
$\frak{I}=\{I \subseteq R|I \textrm{ an ideal of } R\}$ and
$\frak{F}=\{I \subseteq K|I \textrm{ a fractional ideal of } R\}$.
Let $M_\frak{I}=\{f:\frak{I} \rightarrow \frak{I}\}$ and
$M_\frak{F}=\{f:\frak{F} \rightarrow \frak{F}\}$. $M_\frak{I}$  and
$M_\frak{F}$ are clearly a monoids under composition of maps with
identity, the identity map $e:\frak{I} \rightarrow \frak{I}$ and
function composition is associative. $C_R$ will be the subset of
$M_\frak{I}$ consisting of closure operations. Hence, the $f_c$ in
$C_R$ are the set of maps satisfying the following three properties:
(a) $f_c(I) \supseteq I$, (b) $f_c$ preserves inclusions in $R$, and
(c) $f_c\circ f_c=f_c$. $S_R$ will be the set of semiprime
operations of $R$, i.e. $S_R$ are the maps in $C_R$ which also
satisfy (d) $f_c(I)f_c(J) \subseteq f_c(IJ)$. $P_R$ will be the set
of prime operations of $R$, i.e. maps in $S_R$ which also satisfy
(e) $f_c(bI)=bf_c(I)$.

\begin{dfn}  We say a closure operation $f_c$ is bounded on a commutative ring $R$ if for
every maximal ideal $m$ of $R$, there is an $m$-primary ideal $I$ such
that for all $m$-primary $J \subseteq I$, $f_c(J)=I$.  If this is not the
case, we will say that $f_c$ is an unbounded closure operation.
\end{dfn}

The author has determined all the semiprime operations on the set of
ideals of a Dedekind domain and the cuspidal cubic.  In each case,
there are several bounded semiprime operations.  In particular, over
a discrete valuation ring, all the semiprime operations excluding
the identity are bounded.

\begin{prp}[\cite{Va}, Proposition 3.2] When $(R,P)$ is a discrete valuation ring, the
only semiprime operations on the set of ideals of $R$ are the identity,
$f_m$ and $g_m$ defined below:
$$  f_m(P^i)=
\left\{ \begin{array}{l}P^{i} \textrm{ for } 0 \leq i <m\\
             P^m \textrm{ for } i \geq m \end{array}\right.
\textrm{ and } f_m(0)=(0)  $$

and

$$g_m(P^i)=
\left\{ \begin{array}{l}P^{i} \textrm{ for } 0 \leq i <m\\
             P^m \textrm{ for } i \geq m \end{array}\right.
\textrm{ and } g_m(0)=P^m . $$ Moreover, the set of semiprime
operations $S_R$ decomposes as the union $S_R=M_0 \cup M_f$ where
$M_0=\{e\} \cup \{f_n \mid n \geq 0\}$ and $M_f=\{e\} \cup \{g_n
\mid n \geq 0\}$ each of which are submonoids of $M_\frak{I}$ and
$M_f$ is a left but not a right $M_0$-act under composition.
\end{prp}

We make use of the next definition in the following remark.

\begin{dfn} \cite{Va}[Definition 4.10]  Let $R$ be a one-dimensional semigroup ring defined by
$S \subseteq \mathbb{N}_0$.  Let $f_c$ is a bounded semiprime
operation and $J$ be the unique ideal with $f_c(I)=J$ for all $(0)
\neq I  \subseteq J$ and $n \geq 1$ be the conductor of $S$. Suppose
$\mathfrak{a}$ is an ideal which is incomparable to $J$ and
$f_c(\mathfrak{a}) \supseteq J$ and $\mathfrak{a}=\mathfrak{a}_0
\subseteq \mathfrak{a}_1 \subseteq \cdots \subseteq
\mathfrak{a}_k=f_c(\mathfrak{a})$ is a composition series for
$f_c(\mathfrak{a})/\mathfrak{a}$ for $k \geq n$ with $\mathfrak{a}_i
\supseteq J$ for all $i>0$. Then we say $f_c$ is an exceptional
semiprime operation.
\end{dfn}

\begin{rmk}{\rm We would like to note before continuing, that those who study
star operations exclude the zero ideal from the set of fractional
ideals. If we had excluded the zero ideal from the $M_\frak{I}$, the
set of ideals, then we observe that the semiprime operations over a
DVR would be a submonoid of $M_\frak{I}$.  Similarly, if zero were
excluded, the the semiprime operations over a Dedekind domain would
be a submonoid of $M_\frak{I}$ in Proposition 3.6 of \cite{Va}.
However, even excluding zero in the case of the cuspidal cubic,
$S_R$ is not a submonoid because of the exceptional semiprime
operations $f_{n,S,T,m'}$, see Proposition 4.11 in \cite{Va}.}
\end{rmk}

If we consider these same operations on the set of fractional ideals
of $(R,P)$, we see that all the bounded semiprime operations on
$M_{\frak{I}}$ are not semiprime on $M_{\frak{F}}$.  Hence, the only
semiprime operation on the set of fractional ideals of $R$ is the
identity.

\begin{prp}   When $(R,P)$ is a discrete valuation ring, the
only semiprime operation $f \in M_{\frak{F}}$ on the fractional
ideals of $R$ is the identity.
\end{prp}

\begin{proof}
First note, that there are no semiprime operations $f \in
M_{\frak{F}} $ on the set of fractional ideals of $R$, with $R
\subsetneq f(R)$ since if $f(R)=P^{i}$, $i<0$, then
$$f(R)f(R)=P^{2i} \supsetneq P^i=f(R)=f(R \cdot R).$$

The fractional ideals of $R$ have the form $P^i$, $i \in \mathbb{Z}$
and they are totally ordered. Suppose $f_c$ is a bounded semiprime
operation, then $f_c(P^i)=P^m$ for all $i \geq m$.  Suppose $i+j=m$
with $i > m$, then $f_c(P^i)f_c(P^j) \subseteq f_c(P^m)$, but since
$j=m-i<0$, $$P^{m+j} \subseteq f_c(P^i)f_c(P^j) \subseteq
f_c(P^m)=P^m$$ leads to a contradiction since $P^m \subseteq
P^{m+j}$ for $j<0$.
\end{proof}

In fact, no one-dimensional local Noetherian domain can have bounded
semiprime operations on the set of fractional ideals.

\begin{prp}
Let $(R,\frak{m})$ be a local one-dimensional Noetherian domain.
There are no bounded semiprime operations on the fractional ideals
of $R$.
\end{prp}

\begin{proof}
As in the above Proposition, $f_c(R)=R$.  Suppose that $f_c \in M_{\mathfrak{F}}$ is a semiprime operation which is bounded. For any
element $s \in \frak{m}$, and $i \in \mathbb{Z}$, $s^iR$ is a fractional ideal, and the $s^iR$ form a chain.  Suppose  $n$ is the minimal
integer satisfying $s^nR \subseteq I$ for some $m$-primary ideal $I \subseteq R$ and $f_c(s^iR)=I$ for all $i \geq n$ but for $i<n$.  Since
$s^iR \nsubseteq I$, then $f_c(s^iR) \nsubseteq I$. Note that $f_c(s^{i+j}R)=I$ for $i+j=n$.  If $j<0$ and $i>n$, then as $s^{n-j} \notin
I$, then $s^{j}I \supsetneq I$ and the following chain
$$f_c(s^iR)f_c(s^jR)=I f_c(s^jR) \supseteq s^jI \supsetneq I=f_c(s^{i+j}R)$$ implies that $f_c$ is not a semiprime operation since
$f_c(s^jR) f_c(s^iR) \nsubseteq f_c(s^{i+j}R)$.
\end{proof}

Note, there are still semiprime operations over the set of
fractional ideals which are not star operations.  For example, for
any domain which is not Dedekind, the integral closure is a
semiprime operation which is not prime and hence not a star
operation.

\section{Prime operations over one-dimensional semigroup rings}

In \cite{Va}, the author showed that the only prime operation over
the set of ideals of a DVR (Proposition 3.4), a Dedekind domain
(Proposition 3.7), and the cuspidal cubic, $K[[t^2,t^3]]$ (Theorem
4.12) is the identity.  We will show that if $R=K[[t^2,t^{2n+1}]]$,
for $n \geq 1$ then the only prime operation is the identity.
However, we believe that if $R=K[[t^S]]$, for any one-dimensional
semigroup $S$ having ideals minimally generated by three or more
natural numbers, then $R$ has prime operations which are not the
identity.  In fact, we show that there is a non-identity prime
operation on the set of ideals of $K[[t^3,t^4,t^5]]$.

We will determine a nice set of generators for all the ideals in
$R=K[[t^2,t^{2n+1}]]$ and the various relationships between the ideals.
Similar to Proposition 4.1 in \cite{Va}, we will determine a nice form
for the generator of a principal ideal.  Then we will show that
every ideal in $R$ is minimally generated by two elements. First,
we will look at the case $R=K[[t^2,t^5]]$.

\begin{prp}  \label{prid25} Each nonzero nonunit principal ideal of
$R=K[[t^2,t^5]]$ can either be expressed in the form $(t^2+at^5)$, $a \in K$
or $(t^n+at^{n+1}+bt^{n+3})$, $a,b \in K$, $n \geq 4$.
\end{prp}

\begin{proof}
Suppose $0 \neq f \in R$ and $(f) \neq R$.  After multiplying by a
nonzero element of $K$, either
\begin{enumerate}
\item $f=t^2+a_2t^4+a_3t^5+\cdots$ or
\item $f=t^n+a_1t^{n+1}+a_2t^{n+2}+\cdots$ for $n \geq 4$.
\end{enumerate}
Note that in case (1), $n=2$ and $a_1$ must be $0$.  In both cases,
$f=t^2g$ where $g \in K[[t^2,t^3]]$.  We saw in \cite{Va} that if $g \neq
0$ is not a unit, then $(g) = (t^{n-2}+a_1t^{n-1})$ for $n \geq 2$. Thus
$g=(t^{n-2}+at^{n-1})u$, where $u$ is a unit in $K[[t^2,t^3]]$. We would
like to show that $(f)=(t^{n}+at^{n+1}+bt^{n+3})$ for some $a, b \in K$.

First, we will show that $t^{n+2}+a_1t^{n+3} \in (f)$.  Note that
$t^{n+2}+a_1t^{n+3}=t^2(t^{n}+a_1t^{n+1})=(t^{n}+a_1t^{n+1})ut^2u^{-1}=ft^2u^{-1}$.
Now since $u^{-1}$ is a unit in $K[[t^2,t^3]]$, then $t^2u^{-1} \in R$.
Similarly, $t^{n+2}+a_1t^{n+3} \in (t^n+a_1t^{n+1}+(a_3-a_1a_2)t^{n+3})$.

Note if $n \geq 4$, $f=t^nh$ where $h \in K[[t]]$ is a unit.   We will see
that $t^r \in (f)$ for $r \geq n+4$ since $t^r=t^{r-n}h^{-1}f$. Similarly,
$t^r \in (t^n+a_1t^{n+1}+(a_3-a_1a_2)t^{n+3})$.  Since
$f-(t^n+a_1t^{n+1}+(a_3-a_1a_2)t^{n+3}) =
a_2t^{n+2}+a_1a_2t^{n+3}+a_4t^{n+4}+a_5t^{n+5}+\cdots \in (f) \cap
(t^n+a_1t^{n+1}+(a_3-a_1a_2)t^{n+3})$, then
$(f)=(t^n+a_1t^{n+1}+(a_3-a_1a_2)t^{n+3})$.  Hence, either
$(f)=(t^2+at^5)$ for some $a \in K$ or there are $a,b \in K$ with
$(f)=(t^n+at^{n+1}+bt^{n+3})$ for $n \geq 4$.
\end{proof}

Using the fact that the principal ideals have the form $(t^n+at^{n+1}+bt^{n+3})$
we can classify all the non-principal ideals of
$K[[t^2,t^5]]$.

\begin{prp}  \label{ideals25} Each nonzero ideal of $R=K[[t^2,t^5]]$ which is not
principal is minimally generated by two elements of $R$ and can either be expressed in the form $(t^n,t^{n+1})$,  $n \geq 4$, $(t^2,
t^{5})$ , $(t^n+at^{n+1}, t^{n+3})$ or $(t^n,t^{n+1})$ for any $a,b \in K$, $n \geq 4$.
\end{prp}

\begin{proof}
We know the principal ideals are of the form $(t^2+at^5)$ or
$(t^n+at^{n+1}+bt^{n+3})$ by Proposition \ref{prid25}.  Also
$t^{n+r} \in (t^n+at^{n+1}+bt^{n+3})$ for $r \geq 4$, if
$(t^n+at^{n+1}+bt^{n+3}) \in I$ is of minimal initial degree.  Thus,
we need only consider ideals $I$ with additional generators of the
form

\begin{itemize}
\item $t^n+ct^{n+1}+dt^{n+3}$,
\item $t^{n+1}+ct^{n+2}+dt^{n+3}$,
\item $t^{n+2}+ct^{n+3}$ and
\item $ t^{n+3}$.
\end{itemize}

Suppose $I$ has $m-1$ additional generators as above, then as $$V=\{a_0t^n+a_1t^{n+1}+a_2t^{n+2}+a_3t^{n+3} \mid a_i \in K\}$$ is a
$K$-vector subspace of $K[[t^2,t^5]]$, we can put the coefficients in rows of a $m \times 4$ matrix $A$ and put $A$ into reduced row
echelon form to help determine the minimal number of generators of $I$.  The possible reduced row echelon forms for $A$ are

\centerline{(1) $\begin{pmatrix}
1 & a & 0 & b\\
0 & 0 & 0& 0\\
\vdots &\vdots & \vdots & \vdots\\
0 & 0 & 0 &0\\
\end{pmatrix}$, (2) $\begin{pmatrix}
1 & 0 & B & C\\
0 & 1 & D & E\\
0 & 0 & 0 & 0\\
\vdots &\vdots & \vdots & \vdots\\
0 & 0 & 0 & 0\\
\end{pmatrix}$,  (3) $\begin{pmatrix}
1 & B & 0 & C\\
0 & 0 & 1 &D\\
0 & 0 & 0& 0\\
\vdots &\vdots & \vdots& \vdots\\
0 & 0 & 0&0\\
\end{pmatrix}$,
(4) $\begin{pmatrix}
1 & B & C & 0\\
0 & 0 & 0 &1\\
0 & 0 & 0& 0\\
\vdots &\vdots & \vdots& \vdots\\
0 & 0 & 0& 0\\
\end{pmatrix}$,}
\centerline{(5) $\begin{pmatrix}
1 & 0 & 0 & B\\
0 & 1 & 0& C\\
0 & 0 & 1&D\\
0 & 0 & 0& 0\\
\vdots &\vdots & \vdots& \vdots\\
0 & 0 & 0& 0\\
\end{pmatrix}$,
(6) $\begin{pmatrix}
1 & 0 & B & 0\\
0 & 1 & C& 0\\
0 & 0 & 0&1\\
0 & 0 & 0& 0\\
\vdots &\vdots & \vdots& \vdots\\
0 & 0 & 0& 0\\
\end{pmatrix}$,
(7) $\begin{pmatrix}
1 & B & 0 & 0\\
0 & 0 & 1& 0\\
0 & 0 & 0& 1\\
0 & 0 & 0& 0\\
\vdots &\vdots & \vdots& \vdots\\
0 & 0 & 0& 0\\
\end{pmatrix}$,
(8) $\begin{pmatrix}
1 & 0 & 0 & 0\\
0 & 1 & 0& 0\\
0 & 0 & 1& 0\\
0 & 0 & 0& 1\\
0 & 0 & 0& 0\\
\vdots &\vdots & \vdots& \vdots\\
0 & 0 & 0& 0\\
\end{pmatrix}$.}

Now we use the properties of our ring $R$ to simplify further.  (1) represents a principal ideal so we can rule this case out.

(2) represents the ideal $I=(t^n+Bt^{n+2}+Ct^{n+3},t^{n+1}+Dt^{n+2}+Et^{n+3})$.  As
$$t^2(t^n+Bt^{n+2}+Ct^{n+3})-(Bt^{n+4}+Ct^{n+5})=t^{n+2} \in I,$$ then $(t^n+Ct^{n+3},t^{n+1}+Et^{n+3}) \subseteq I$.
Similarly,
$$t^2(t^{n+1}+Et^{n+3})-Et^{n+5}=t^{n+3} \in I$$ implies that $(t^n,t^{n+1}) \subseteq I$.
Clearly, $t^{n+2},t^{n+3} \in (t^n,t^{n+1})$ also implying that
$I=(t^n,t^{n+1})$.

(3) represents the ideal $I=(t^n+Bt^{n+1}+Ct^{n+3},t^{n+2}+Dt^{n+3})$. As
$$t^2(t^n+Bt^{n+1}+Ct^{n+3})-(t^{n+2}+Dt^{n+3})-Ct^{n+5}=(B-D)t^{n+3} \in I,$$ we now observe that if $B-D \neq 0$, then $t^{n+3}$ and hence
$t^{n+2}$ and $t^n+Bt^{n+1}$ are in $I$ and clearly $I \subseteq (t^n+Bt^{n+1},t^{n+3})$ implying that $I=(t^n+Bt^{n+1},t^{n+3})$.  In the
case that $B-D=0$, then $I$ is principal generated by $t^n+Bt^{n+1}+Ct^{n+3}$.

(4) represents the ideal $I=(t^n+Bt^{n+1}+Ct^{n+2},t^{n+3})$.  Since
$$t^{2}(t^{n}+Bt^{n+1}+Ct^{n+3})-(Bt^{n+3}+Ct^{n+5})=t^{n+2} \in I$$
and similarly $t^{n+2} \in (t^n+Bt^{n+1}+Ct^{n+2},t^{n+3})$ using
similar reasoning to (2) above we see that
$I=(t^n+Bt^{n+1},t^{n+3})$.

(5) represents the ideal $I=(t^n+Bt^{n+3},t^{n+1}+Ct^{n+3},t^{n+2}+Dt^{n+3})$.  Note that
$$t^{n+2}+Dt^{n+3}=t^2(t^{n}+Bt^{n+3})+Dt^2(t^{n+1}+Ct^{n+3})-(B+DC)t^{n+5}.$$ As we argued in (2) above, we
see that $I=(t^n,t^{n+1})$.

(6) represents the ideal $I=(t^n+Bt^{n+2},t^{n+1}+Ct^{n+2},t^{n+3})$.  As $$t^{n+2}=t^2(t^n+Bt^{n+2})-Bt^{n+4} \in I$$ and
$t^{n+3}=t^2{t^{n+1}}$, $I$ is minimally generated by $(t^n,t^{n+1})$.

(7) represents the ideal $I=(t^n+Bt^{n+1},t^{n+2},t^{n+3})$.  Since $t^{n+2}=t^2(t^n+Bt^{n+2})-Bt^{n+3}$, we see that
$I=(t^n+Bt^{n+1},t^{n+3})$.

(8) represents the ideal $I=(t^n,t^{n+1},t^{n+2},t^{n+3})$. Clearly,
$I$ is minimally generated by $t^n$ and $t^{n+1}$.
\end{proof}

In fact a lattice of the ideals of $K[[t^2,t^5]]$ is as follows:
$$ 
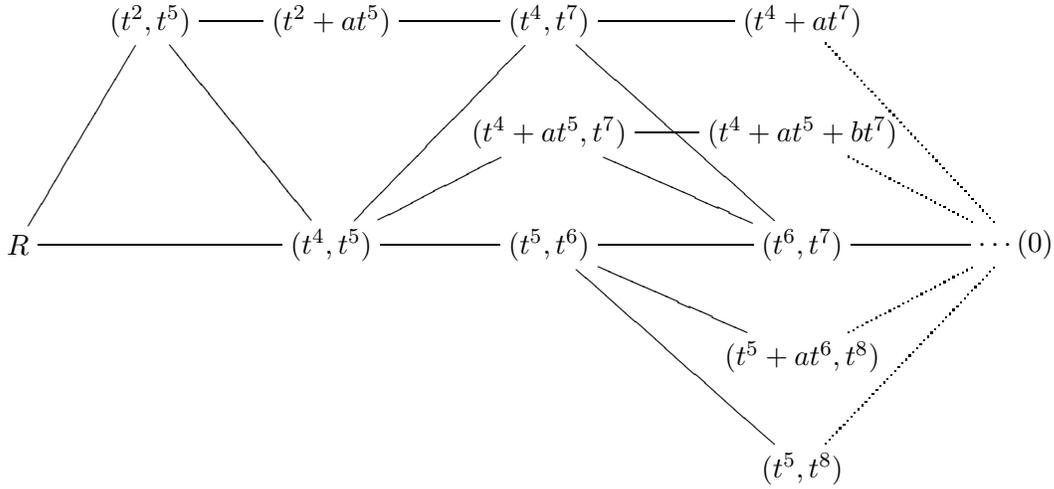
\xymatrix{ & (t^2,t^5)\ar @{-} [ddr]  \ar @{-} [r]    &   (t^2+at^5)  \ar @{-} [r]      & (t^4,t^7) \ar @{-} [ddl] \ar @{-} [ddr]\ar @{-} [r]&  (t^4+at^7) \ar @{.} [ddr]&  & \\
 &     &     & (t^4+at^5,t^7) \ar @{-} [dl] \ar @{-} [dr]\ar @{-} [r]&  (t^4+at^5+bt^7)\ar @{.} [dr]&  & \\
    R\ar @{-} [rr] \ar @{-} [uur]& & (t^4,t^5)\ar @{-} [r]& (t^5,t^6)\ar @{-} [r] \ar @{-} [dr]&(t^6,t^7)\ar @{-} [r]&\cdots  (0)\ar  @{.} [l]\\
                 &       &    &      &   (t^5+at^6,t^8)   \ar @{.} [ur]&\\
  &       &    &      &   (t^5,t^8)   \ar @{-} [uul] \ar @{.} [uur]&\\}$$ where the lines indicate $\supseteq$. At each node where
there is an $a \in K$ present in the expression of the ideal, there are a cardinality of $K$'s worth of ideals at that node.  If both $a,b
\in K$ are present, there are a cardinality of $K^2$'s worth of ideals at that node.

We can similarly determine the ideal structure of $K[[t^2,t^{2r+1}]]$.  First we express the principal ideals similar to the way we
expressed them in Proposition \ref{prid25}.  In fact, the proof is very similar, building up on rings of the form $K[[t^2,t^{2m+1}]]$ for
$m<r$.

\begin{prp}  \label{prid2odd} Each nonzero nonunit principal ideal of $R=K[[t^2,t^{2r+1}]]$ can be expressed
in the form $(t^2+a_1t^{2r+1}), (t^4+a_1t^{2r+1}+a_3t^{2r+3}), \ldots, (t^{2r-2}+a_1t^{2r+1}+\cdots a_{2r-3}t^{4r-3})$,  or
$(t^n+a_1t^{n+1}+a_3t^{n+3} + \cdots+ a_{2r-1}t^{n+2r-1})$, $n \geq 2r$ and $a_{i} \in K$.
\end{prp}

\begin{proof}
The proof is similar to Proposition \ref{prid25}. Suppose $0 \neq f \in R$ and $(f) \neq R$.  After multiplying by a nonzero element of
$K$, either
\begin{enumerate}
\item $f=t^{2k}+a_2t^{2k+2}+\cdots+a_{2(r-k)}t^{2r}+a_{2(r-k)+1}t^{2r+1}+\cdots$ for $1 \leq k<r$ or
\item $f=t^n+a_1t^{n+1}+a_2t^{n+2}+\cdots$ for $n \geq 2r$.
\end{enumerate}
Note that in case (1), $n=2k$ for $1 \leq k<r$ and $a_1=a_3=a_{2(r-k)+1}=0$. In both cases, $f=t^2g$ where $g \in K[[t^2,t^{2r-1}]]$.  We
assume by induction that the principal ideals of $K[[t^2,t^{2r-1}]]$ are of the form $(t^2+a_1t^{2r-1}), (t^4+a_1t^{2r-1}+a_3t^{2r+1}),
\ldots, (t^{2r-4}+a_1t^{2r-1}+\cdots a_{2r-5}t^{4r-7})$,  or $(t^n+a_1t^{n+1}+a_2t^{n+3} + \cdots+ a_{2r-3}t^{n+2r-3})$, $n \geq 2(r-1)$
and $a_{i} \in K$. As in Proposition \ref{prid25} we see that if $g \neq 0$ is not a unit, then $(g) = (t^{n-2}+a_{1}1t^{n-1}+\cdots
a_{2r-5}t^{n+2r-3})$ for $n \geq 2$. Thus $g=(t^{n-2}+a_1t^{n-1}+a_3t^{n+1}+\cdots a_{2r-5}t^{n+2r-3} )u$, where $u$ is a unit in
$K[[t^2,t^{2r-1}]]$. We would like to show that $(f)=(t^{n}+a_1t^{n+1}+a_3t^{n+3}+ \cdots a_{2r-1}t^{n+2r-1})$ for $a_i\in K$.

First, we will show that $t^{n+2}+a_1t^{n+3}+a_{2r-3}t^{n+2r-1} \in
(f)$.  Note that
\begin{eqnarray*}ft^2u^{-1}&=&t^2(t^{n}+a_1t^{n+1}+a_{2r-3}t^{n+2r-3})ut^2u^{-1}\\
&=&t^{n+2}+a_1t^{n+3}+a_{2r-3}t^{n+2r-1}.\end{eqnarray*} Now since
$u^{-1}$ is a unit in $K[[t^2,t^{2r-1}]]$, then $t^2u^{-1} \in R$.
Similarly, $$t^{n+2}+a_1t^{n+3}+a_{2r-3}t^{n+2r-1} \in
(t^n+a_1t^{n+1}+b_3t^{n+3}+\cdots b_{2r-1}t^{2r-1})$$ where
$b_3=a_3-a_2a_1$, $b_5=a_5-a_4a_1-a_2a_3+a_1(a_2)^2$ and generally
$$b_{2m-1}=a_{2m-1}-\sum\limits_{i=1}^{m-1}
a_{2i-1}a_{2(m-i)}+\sum\limits_{i=1}^{m-2} \sum\limits_{j=1}^{m-i}
a_{2i-1}a_{2j}a_{2(m-i-j)}-\cdots +(-1)^{m-1}a_1a_2^{m-1}.$$

Note if $n \geq 2r$, $f=t^nh$ where $h \in K[[t]]$ is a unit.   We will see that $t^m \in (f)$ for $m \geq n+2r$ since
$t^m=t^{m-n}h^{-1}f$. Similarly, $t^m \in (t^n+a_1t^{n+1}+(a_3-a_1a_2)t^{n+3})$.  Since $f-(t^n+a_1t^{n+1}+b_3t^{n+3}+\cdots +
b_{2r-1}t^{n+2r-1})  \in (f) \cap (t^n+a_1t^{n+1}+b_3t^{n+3}+ \cdots +b_{2r-1}t^{n+2r-1})$, then $(f)=(t^n+a_1t^{n+1}+b_3t^{n+3}+ \cdots
+b_{2r-1}t^{n+2r-1})$. Hence, the principal ideals are those stated in the Proposition.
\end{proof}

Using the structure of the principal ideals which we determined in Proposition \ref{prid2odd}, we can now determine the structure of all
non-principal ideals.

\begin{prp}  \label{ideals2odd} Each nonzero ideal of $R=K[[t^2,t^{2r+1}]]$ which is not principal is minimally generated
by two elements of $R$ and can either be expressed in the form $(t^2,t^{2r+1})$, $(t^4,t^{2r+1})$, $(t^4+a_1t^{2r+1},t^{2r+3}), \ldots,
(t^{2r-2}+a_1t^{2r+1}+\cdots+a_{r-1}t^{4r-5}, a_{r}t^{4r-3})$,  or $(t^n, t^{n+1})$, $((t^n+a_1t^{n+1},t^{n+3})$,
$(t^n+a_1t^{n+1}+a_3t^{n+3} + \cdots+ a_{2r-3}t^{n+2r-3},t^{n+2r-1})$, $n \geq 2r$ and $a_{2m+1} \in K$ for $0 \leq m \leq r$.
\end{prp}

\begin{proof}
The proof is similar to the proof of Proposition \ref{ideals25},
albeit slightly more technical. We will only give a sketch of it
here.  As in the proof of Proposition \ref{ideals25}, we can
determine the possible generators for an ideal $I$ containing
$f=t^n+a_1t^{n+1}+\cdots+a_rt^{n+2r-1}$ where $f$ was an element of
$I$ with a minimal initial term.  Then we will put the coefficients
in an $m \times 2r$ matrix $A$ and determine all reduced row echelon
expressions for $A$.  After we do this, we use the properties of the
ring to see that all the ideals will have the structure as in the
statement of the Proposition.
\end{proof}

Now we will see that the only prime operation on the ring $K[[t^2,t^{2r+1}]]$ is the identity.

\begin{prp}
The only prime operation on the set of ideals of
$R=K[[t^2,t^{2r+1}]]$ is the identity.
\end{prp}

\begin{proof}
Since every prime operation satisfies \begin{equation}\label{primeeq} f_c((g)I)=(g)f_c(I)\end{equation} for all principal ideals $(g) \in
R$, if we substitute $R$ in for $I$ in (\ref{primeeq}), we obtain $f_c((g))=(g)$.  Thus all principal ideals must be $f_c$-closed.  Note
that each principal ideal is sandwiched between many 2-generated ideals. Suppose $f_c(I) \neq I$ for some two generated ideal
$$I=(t^n+a_1t^{n+1}+\cdots + a_{2m-3}t^{n+2m-3},t^{n+2m-1})$$ for $2\leq m \leq r$.  Then $f_c((g)I)= (g)f_c(I) \neq (g)I$ for any principal
ideal $(g)$.  Since the conductor of $R$ is $\mathfrak{c}=(t^{2r},
\ldots, t^{4r-1})$.  Then for any $(g) \in \mathfrak{c}$, $f_c((g)I)
\neq (g)I$.  Hence, all the ideals of the form
$J_{k,S}=(t^k+b_1t^{k+1}+\cdots + b_{2m-3}t^{k+2m-3},t^{k+2m-1})$
for $k \geq n+2r$ and $S$ the $m-1$-tuple $( b_1, \ldots, b_{m-1})$
satisfy $f_c(J_{k,S}) \neq J_{k,S}$. Note $$f_c(J_{k,S})\supseteq
(t^{k-2r+2m-2}+b_1t^{k-2r+2m-1}+\cdots +
b_{2m-3}t^{k-2r+4m-5}+bt^{k-2r+4m-3})$$ some $b \in K$ or
$$f_c(J_{k,S})\supseteq (t^k+b_1t^{k+1}+\cdots +
b_{k+2m-5}t^{k+2m-5},t^{k+2m-3}).$$

If $f_c(J_{k,S})=(t^k+b_1t^{k+1}+\cdots +
b_{k+2m-5}t^{k+2m-5},t^{k+2m-3})$, then since
$$(t^{k-2r+2m-2}+b_1t^{k-2r+2m-1}+\cdots +
b_{2m-3}t^{k-2r+4m-5}+bt^{k-2r+4m-3})$$ is incomparable to
$$(t^k+b_1t^{k+1}+\cdots + b_{k+2m-5}t^{k+2m-5},t^{k+2m-3})$$ and
$$f_c((t^{k-2r+2m-2}+b_1t^{k-2r+2m-1}+\cdots +
b_{2m-3}t^{k-2r+4m-5}+bt^{k-2r+4m-3}))$$ contains
$(t^k+b_1t^{k+1}+\cdots + b_{k+2m-5}t^{k+2m-5},t^{k+2m-3}),$ this
forces
$$f_c((t^{k-2r+2m-2}+b_1t^{k-2r+2m-1}+\cdots +
b_{2m-3}t^{k-2r+4m-5}+bt^{k-2r+4m-3}))$$ to be not equal to
$$(t^{k-2r+2m-2}+b_1t^{k-2r+2m-1}+\cdots +
b_{2m-3}t^{k-2r+4m-5}+bt^{k-2r+4m-3}),$$ contradicting that
$f_c(g)=(g)$ for all principal $(g)$.

If $f_c(J_{k,S})=(t^{k-2r+2m-2}+b_1t^{k-2r+2m-1}+\cdots +
b_{2m-3}t^{k-2r+4m-5}+bt^{k-2r+4m-3})$ then
$$f_c((t^k+b_1t^{k+1}+\cdots + b_{k+2m-5}t^{k+2m-5},t^{k+2m-3}))$$
contains $$(t^{k-2r+2m-2}+b_1t^{k-2r+2m-1}+\cdots +
b_{2m-3}t^{k-2r+4m-5}+bt^{k-2r+4m-3})$$ which will also force
$$f_c((t^{k+2r-2m+4}+b_1t^{k+2r-2m+5}+\cdots +
b_{k+2m-5}t^{k+2r-1},t^{k+2r+1}))$$ to contain
$$(t^{k+2}+b_1t^{k+3}+\cdots + b_{2m-3}t^{k+2m-1}+bt^{k+2m+1}).$$ Now
since $(t^k+b_1t^{k+1}+\cdots + b_{k+2m-5}t^{k+2m-5},t^{k+2m-3})$ is
incomparable to $$(t^{k-2r+2m-2}+b_1t^{k-2r+2m-1}+\cdots +
b_{2m-3}t^{k-2r+4m-5}+bt^{k-2r+4m-3})$$ and the fact that
$(t^k+b_1t^{k+1}+\cdots + b_{k+2m-5}t^{k+2m-5},t^{k+2m-3})$ contains
 $$(t^k+b_1t^{k+1}+\cdots + b_{k+2m-5}t^{k+2m-5}+ct^{k+2m-3})$$ which
contains $$ (t^{k+2r-2m+4}+b_1t^{k+2r-2m+5}+\cdots +
b_{k+2m-5}t^{k+2r-1},t^{k+2r+1}),$$ then we see that
$f_c((t^k+b_1t^{k+1}+\cdots + b_{k+2m-5}t^{k+2m-5}+ct^{k+2m-3}))$ is
not equal to $$(t^k+b_1t^{k+1}+\cdots +
b_{k+2m-5}t^{k+2m-5}+ct^{k+2m-3})$$ some $c \neq b_{2m-3}$, again
contradicting that $f_c(g)=(g)$ for all principal $(g)$.
\end{proof}

Although, there are no non-identity prime operations over semigroup
rings which have ideals minimally generated by at most two elements.
We can show that there is a prime operation on $k[[t^3,t^4,t^5]]$
which is not the identity. This answers a question posed on the last
page of \cite{Va} in the negative. We suspect that in semigroup
rings which have ideals minimally generated by 3 or more elements
will also have prime operations which are not the identity. First,
we generators for all the ideals of $k[[t^3,t^4,t^5]]$.

\begin{prp}\label{ideals3}
Let $R=k[[t^3,t^4,t^5]]$.  Then the nonzero, nonunit ideals of $R$ can be expressed in the forms $(t^n+at^{n+1}+bt^{n+2}),
(t^n+at^{n+1},t^{n+2}), (t^n+at^{n+2},t^{n+1}+bt^{n+2})$  or $(t^n,t^{n+1},t^{n+2})$ for $n\geq 3$ and $a,b \in k$.
\end{prp}

\begin{proof}
In \cite{FV}, Fouli and the current author showed in Proposition 5.2 that all nonunit principal ideals are of the form
$(t^n+at^{n+1}+bt^{n+2})$, for $n\geq 3$ and $a,b \in k$.

Now we consider the form of all non-principal ideals $I$.  Since $(t^{n+3},t^{n+4},t^{n+5}) \subseteq (t^n+at^{n+1}+bt^{n+2})$, if
$(t^n+at^{n+1}+bt^{n+2}) \in I$ is the form with the initial term of smallest degree, the other generators can only be of the form
$t^n+ct^{n+1}+dt^{n+2}$, $t^{n+1}+et^{n+2})$ and $t^{n+2})$.

Observe that $V=\{a_0t^n+a_1t^{n+1}+a_2t^{n+2} \mid a_i \in k, n \geq 3\}$ is a subspace of the $k$-vector space $k[[t^3,t^4,t^5]]$.  To
determine a minimal set of generators for $I$, we can set the coefficients of the $m$ generators of $I$ in rows of an $m \times 3$ matrix
$A$.  The possible reduced row echelon forms for $A$ are:

\centerline{$\begin{pmatrix}
1 & a & b\\
0 & 0 & 0\\
\vdots &\vdots & \vdots\\
0 & 0 & 0\\
\end{pmatrix}$, $\begin{pmatrix}
1 & 0 & B\\
0 & 1 & C\\
0 & 0 & 0\\
\vdots &\vdots & \vdots\\
0 & 0 & 0\\
\end{pmatrix}$,  $\begin{pmatrix}
1 & B & 0\\
0 & 0 & 1\\
0 & 0 & 0\\
\vdots &\vdots & \vdots\\
0 & 0 & 0\\
\end{pmatrix}$, $\begin{pmatrix}
1 & 0 & 0\\
0 & 1 & 0\\
0 & 0 & 1\\
0 & 0 & 0\\
\vdots &\vdots & \vdots\\
0 & 0 & 0\\
\end{pmatrix}$}

\noindent for $B$ and $C$ any elements of $k$.  These row echelon forms correspond to the ideals listed in the statement of the
proposition.
\end{proof}

\begin{prp}
There is a prime operation $f_c$ on the set of ideals of $R=k[[t^3,t^4,t^5]]$ which is not the identity.
\end{prp}

\begin{proof}
The principal ideal $(f)=(t^n+at^{n+1}+bt^{n+2})$ is contained in all the two-generated ideals containing
$$(t^n,t^{n+1}+\frac{b}{a}t^{n+2}), (t^n+at^{n+1},t^{n+2}), (t^n+ct^{n+2},t^{n+1}+\frac{b-c}{a}t^{n+2}),$$ for all $c \in k$ as long as $0
\neq a \in k$.  $(f)$ is also contained in all the three generated ideals containing $(t^{n},t^{n+1},t^{n+2})$.  The ideal
$(t^n+at^{n+1}+bt^{n+2})$ contains all ideals contained in $(t^{n+3},t^{n+4},t^{n+5})$.

The following ideals are incomparable to $(t^n+at^{n+1}+bt^{n+2})$: \begin{itemize}
\item 3-generated ideals: $(t^{n+1},t^{n+2},t^{n+3})$,
$(t^{n+2},t^{n+3},t^{n+4})$,

\item 2-generated ideals with initial degree greater than $n$: $(t^{n+1},t^{n+2}+dt^{n+3})$,
$(t^{n+1}+dt^{n+2},t^{n+3})$, $(t^{n+1}+dt^{n+3},t^{n+2}+et^{n+3})$,
$(t^{n+2},t^{n+3}+dt^{n+4})$, $(t^{n+2}+dt^{n+3},t^{n+4})$ and $(t^{n+2}+dt^{n+4},t^{n+3}+et^{n+4})$  for any $d, e \in k$,

\item 2-generated ideals with initial degree less than or equal to $n$:  $(t^n,t^{n+1}+dt^{n+2})$ or $(t^{n-2}+dt^{n-1},t^{n})$ with
$d \neq \frac{b}{a}$, $(t^n+dt^{n+1},t^{n+2})$ with $d \neq a$, $(t^n+dt^{n+2},t^{n+1}+et^{n+2})$ with $ae \neq b-d$,
$(t^{n-1}+dt^n,t^{n+1})$ for any $d \in k$, $(t^{n-1}+dt^{n+1},t^n+et^{n+2})$ for any $d \in k$ and $a \neq e \in k$ and
$(t^{n-2}+dt^{n},t^{n-1}+et^n)$ for any $d,e \in k$.

\item principal ideals: $(t^n+dt^{n+1}+et^{n+2})$ for $(a,b) \neq(d,e) \in k^2$, $(t^{n+i}+dt^{n+i+1}+et^{n+i+2})$ for
any $d,e \in k$ and $i \in \{-2,-1,1,2\}$.

\end{itemize}

Depicting a lattice of all ideals and their relationships with each other would be quite complicated, but we can easily draw the lattice of
three generated ideals and principal ideals:

$$ 
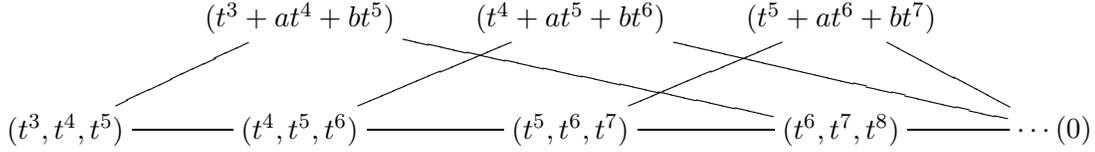
\xymatrix{
                         &   (t^3+at^4+bt^5)  \ar @{-} [dl] \ar@{-}[drr] & (t^4+at^5+bt^6)  \ar @{-} [dl] \ar@{-}[drr]&  (t^5+at^6+bt^7)  \ar @{-} [dl] \ar@{-}[dr]&\\
(t^3,t^4,t^5)\ar @{-} [r]& (t^4,t^5,t^6)\ar @{-} [r]   &(t^5,t^6,t^7)\ar @{-} [r]       &(t^6,t^7,t^8)\ar @{-} [r]              &\cdots
(0)\ar  @{.} [l]\\}$$ where each of the lines denote $\supseteq$ and at each principal node, there are the cardinality of $k^2$ ideals.

We will use this lattice to define a prime operation $f_c$ on the ideals of $R$ as follows: $f_c(I)=(t^n,t^{n+1},t^{n+2})$ if
$I=(t^n,t^{n+1}+ct^{n+2})$, $I=(t^n+dt^{n+1},t^{n+2})$, $I=(t^n+dt^{n+1},t^{n+1}+et^{n+2})$ or $I=(t^n,t^{n+1},t^{n+2})$.  However,
$f_c((t^n+at^{n+1}+bt^{n+2}))=(t^n+at^{n+1}+bt^{n+2})$.

Now we verify that $f_c$ is a prime operation.  Clearly, $I \subseteq f_c(I)$.  Also since any ideal $I$ containing
$(t^n+at^{n+1}+bt^{n+2})$ has the property that $f_c(I) \supseteq (t^n,t^{n+1},t^{n+2})$ and any ideal contained in
$(t^n+at^{n+1}+bt^{n+2}) \supseteq (t^{n+3},t^{n+4},t^{n+5})$ has the property that $f_c(I) \subseteq (t^{n+3},t^{n+4},t^{n+5}) \subseteq
(t^n+at^{n+1}+bt^{n+2})$ then if $I \subseteq J$, then $f_c(I) \subseteq f_c(J)$.  Clearly, $f_c(f_c(I))=f_c(I)$ because if $I$ is
principal then $f_c(I)=I$ and if $I$ is not principal then $f_c(I)=(t^n,t^{n+1},t^{n+2})$ for some $n$ and $(t^n,t^{n+1},t^{n+2})$ is
$f_c$-closed.

Note that the product of two principal ideals is principal, thus $f_c((f)(g))=f_c((f))f_c((g))$ for principal ideals $(f)$ and $(g)$.  The
product of any ideal with a non-principal ideal will be non-principal.  If $I=(t^n+at^{n+1}+bt^{n+2})+I^{\prime}$ and
$J=(t^m+ct^{m+1}+dt^{m+2})+J^{\prime}$ with $I^{\prime}$ incomparable to $(t^n+at^{n+1}+bt^{n+2})$ and $J^{\prime}$ incomparable to
$(t^m+ct^{m+1}+dt^{m+2})$ or the zero ideal, then $IJ=(t^{m+n}+(a+c)t^{m+n+1}+(ac+b+d)t^{m+n+2})+K^{\prime}$ where
$K^{\prime}=(t^n+at^{n+1}+bt^{n+2})J^{\prime}+(t^m+ct^{m+1}+dt^{m+2})I^{\prime}+I^{\prime}J^{\prime}$ which is minimally generated by
either 2 or 3 elements.  In any case $f_c(IJ)=(t^{m+n},t^{m+n+1},t^{m+n+2})=f_c(I)f_c(J)$.  Hence, $f_c$ is a prime operation which is not
the identity.

\end{proof}

It seems very likely, that we should be able to exhibit non-identity
prime operations in one-dimensional semigroup rings with ideals
minimally generated by three or more elements using an argument
similar to the one above.

\end{document}